\documentclass[11pt,  reqno]
{amsart}

\usepackage{amsmath,amssymb,amscd,amsthm,amsxtra, esint}

\allowdisplaybreaks[2]

\newtheorem{theorem}{Theorem} [section]

\newtheorem{oldtheorem}{Theorem}

\newtheorem{maintheorem}{Theorem}

\newtheorem{proposition}[theorem]{Proposition}
\newtheorem{remark}[theorem]{Remark}

\newtheorem{definition}[theorem]{Definition}

\DeclareMathOperator*{\supp}{supp}

\DeclareMathOperator*{\ave}{ave}

\DeclareMathOperator{\diam}{diam}

\newcommand{\I}{\hspace{0.5mm}\text{I}\hspace{0.5mm}}
\newcommand{\II}{\text{I \hspace{-2.8mm} I} }
\newcommand{\III}{\text{I \hspace{-2.9mm} I \hspace{-2.9mm} I}}

\newcommand{\IV}{\text{I \hspace{-2.9mm} V}}

\newcommand{\noi}{\noindent}
\newcommand{\Z}{\mathbb{Z}}
\newcommand{\R}{\mathbb{R}}
\newcommand{\C}{\mathbb{C}}

\let\Re=\undefined\DeclareMathOperator*{\Re}{Re}

\newcommand{\al}{\alpha}
\newcommand{\be}{\beta}
\newcommand{\dl}{\delta}

\newcommand{\Dl}{\Delta}
\newcommand{\eps}{\varepsilon}

\newcommand{\wt}{\widetilde}

\newcommand{\dd}{\partial}

\renewcommand{\l}{\ell}

\newcommand{\les}{\lesssim}

\newcommand{\jb}[1]
{\langle #1 \rangle}

\newcommand{\ind}{\mathbf 1}

\renewcommand{\S}{\mathcal{S}}

\numberwithin{equation}{section}
\numberwithin{theorem}{section}

\newcommand{\BMO}{\textit{BMO}}
\newcommand{\CZ}{Calder\'on-Zygmund }

\begin{document}

\baselineskip = 14pt


\title[$T(b)$ theorems \`a la Stein]
{Linear and bilinear $T(b)$ theorems
\`a la Stein}

\author{\'Arp\'ad  B\'enyi and Tadahiro Oh}

\address{
\'Arp\'ad  B\'enyi\\
Department of Mathematics\\
 Western Washington University\\
 516 High Street, Bellingham\\
  WA 98225\\ USA}
\email{arpad.benyi@wwu.edu}

\address{
Tadahiro Oh\\
School of Mathematics\\
The University of Edinburgh,
and The Maxwell Institute for the Mathematical Sciences\\
James Clerk Maxwell Building\\
The King's Buildings\\
Peter Guthrie Tait Road\\
Edinburgh\\ EH9 3FD\\ United Kingdom}
\email{hiro.oh@ed.ac.uk}

\begin{abstract}
In this work, we state and prove versions of
the linear and bilinear $T(b)$ theorems
involving quantitative estimates, analogous
to the quantitative linear $T(1)$ theorem due to Stein.

\end{abstract}



\maketitle

\vspace{-5mm}

\section{Introduction}

The impact of the classical Calder\'on-Zygmund theory permeates through analysis and PDEs. Nowadays, both the linear and multilinear aspects of this theory are well understood and continue to be intertwined with aspects of analysis that are beyond their reach, such as those considering the bilinear Hilbert transform.

Two fundamental results in the linear theory from the 1980's are the celebrated $T(1)$ theorem of David and Journ\'e \cite{DJ} and $T(b)$ theorem of David, Journ\'e, and Semmes \cite{DJS}.
Both results were 
strongly motivated by
the study of
the Cauchy integral on a Lipschitz curve and the related Calder\'on commutators. Their gist lies in understanding the boundedness of a singular operator via appropriate simpler testing conditions.

In the $T(1)$ theorem, one needs to test a singular operator and its transpose on the constant function 1. If both the operator and its transpose were $L^\infty\to \BMO$ bounded, then by duality and interpolation \cite{FS}, the operator would be bounded on $L^2$. The remarkable aspect of the $T(1)$ theorem is that one does not need to test the operator on the whole $L^\infty$, but just on one special element in it. Going back to the Cauchy integral operator associated to a Lipschitz function $A$, it turns out that it is not necessarily easy to test the operator on 1. It is, however, much easier to test the operator on the $L^\infty$ function $1+iA'$. Thus, as the name suggests, the $T(b)$ theorem extends the $T(1)$ theorem by replacing the constant function 1 with a suitable $L^\infty$ function $b$; or, to be more precise, by replacing 1 with two suitable functions $b_0$ and $b_1$ in $L^\infty$ on which we test an operator and its transpose. The bilinear Calder\'on-Zygmund theory has its own versions of the $T(1)$ and $T(b)$ theorems, such as those proved by Grafakos and Torres \cite{GT1} and by  Hart \cite{Hart2}, respectively.
See Theorems \ref{THM:BT1} and \ref{THM:BTB} below.

In  this work, we revisit the $T(b)$ theorem, both in linear and bilinear setting, through the lens of a gem due to Stein \cite{Stein}. We are alluding to his formulation of the $T(1)$ theorem involving quantitative estimates for a singular  operator and its transpose when tested now on normalized bump functions.
Our goal is to prove that an analogous  natural formulation \`a la Stein can be given
for  the $T(b)$ theorems in the linear and bilinear settings. We note that, while for the sake of clarity in our presentation we have chosen to delineate the linear and bilinear settings, a unified discussion is certainly possible under the encompassing more general multilinear setting.

\section{Linear \CZ theory}
\label{SEC:linear}

In this section,  we consider a linear  singular operator $T$
a priori defined from $\S$ into $\S'$
of the form
\begin{align}
T(f)(x) = \int_{\R^d} K(x, y) f(y) dy .
\label{CZ1}
\end{align}

\noi
Here, we assume that, away from the diagonal
 $\Dl = \{ (x, y) \in \R^{2d}:\, x = y\}$,
 the distributional kernel $K$ of $T$ coincides with
a function that is locally integrable on $\R^{2d} \setminus \Dl$.
The formal transpose $T^*$ of $T$
is defined similarly
with
the kernel  $K^*$
given by
$K^*(x, y) : = K(y, x)$.

\begin{definition}
\rm

A locally integrable function $K$ on $ \R^{2d} \setminus \Dl$
is called a (linear) {\it \CZ  kernel}
if it satisfies the following conditions.

\begin{itemize}
\item[\textup{(i)}]
For all $x, y \in \R^d$, we have
$|K(x, y)|  \les |x-y|^{-d}$,

\smallskip

\item[\textup{(ii)}]
There exists $\dl \in (0, 1]$ such that
\begin{align}
|K(x, y) - K(x', y)|
& + |K(y, x) - K(y, x')|
\les \frac{|x - x'|^\dl}{|x-y|^{d+\dl}}
\label{CZ2}
\end{align}

\noi
for all $x, x', y \in \R^d$
satisfying
$|x - x'| < \frac 12 |x-y|$.

\end{itemize}

\end{definition}

\noi
We say that a linear singular operator $T$
of the form \eqref{CZ1}
with a \CZ kernel
is a {\it linear \CZ operator} if $T$ extends
to a bounded operator on $L^{p_0}$
for some $1< p_0 < \infty$.
It is well known \cite{Stein1} that
if  $T$ is a linear \CZ operator,
then it is bounded on $L^p$ for all
$1< p < \infty$. Hence, in the following, we restrict our attention to the $L^2$-boundedness of such linear operators. We point out that the \CZ operator $T$ is also $L^\infty\to \BMO$ bounded.
Here, $\BMO$ denotes the space of functions of \emph{bounded mean oscillation}, which we now recall.

\begin{definition} \rm
Given a locally integrable function $f$ on $\R^d$,
define the $\BMO$-seminorm by
\[ \|f\|_{\BMO} := \sup_{Q} \frac{1}{|Q|} \int_Q|f(x) - \ave_Q f|dx,\]
	
\noi
where the supremum is taken over all cubes $Q \subset \R^d$
and \[\ave_Q f : = \frac{1}{|Q|} \int_Q f(x) dx.\]
Then, we say that $f$ is of bounded mean oscillation if $\|f\|_{\BMO}< \infty$
and
we define $\BMO(\R^d)$ by
\[ \BMO(\R^d) : =
\big\{ f \in L^1_{\text{loc}}(\R^d):\,
\| f \|_{\BMO} < \infty \big\}.\]

\end{definition}

\subsection{Classical linear $T(1)$ and $T(b)$ theorems}

In this subsection, we provide a brief discussion of the classical $T(1)$ and $T(b)$ theorems proved in \cite{DJ} and \cite{DJS}, respectively. In order to do so, we need to define a few more notions.

\begin{definition}\rm
We say that a function $\phi \in \mathcal{D}$ is
a \emph{normalized bump function} of order $M$ if
$\supp \phi \subset B_0(1)$
and $\|\dd^\al \phi\|_{L^\infty} \leq 1$
for all multi-indices $\al$ with $|\al| \leq M$.

\end{definition}

\noi
Here, $B_x(r)$ denotes the ball of radius $r$ centered at $x$.
Given $x_0 \in \R^d$ and $R>0$, we set
\begin{align}
\phi^{x_0, R}(x) = \phi\Big(\frac{x-x_0}{R}\Big).
\label{scaling}
\end{align}

\begin{definition} \label{DEF:WBP1} \rm

We say that a linear singular integral operator $T: \S\to \S'$
has the {\it weak boundedness property}
if there exists $M \in \mathbb N\cup \{0\}$
such that we have
\begin{align}
 \big| \big\langle T(\phi_1^{x_1, R}), \phi_2^{x_2, R}\big\rangle \big|\les R^d
\label{WBP1}
 \end{align}

\noi
for all normalized bump functions $\phi_1$ and $\phi_2$
of order $M$, $x_1, x_2 \in \R^d$,
and $R>0$.

\end{definition}

\noi
We note that it suffices to verify \eqref{WBP1} for $x_1 = x_2$;
see \cite{Hart1}. The statement of the $T(1)$ theorem of David and Journ\'e \cite{DJ} is the following.

\begin{oldtheorem}[$T(1)$ theorem]\label{THM:DJ}
Let $T: \S \to \S'$ be a linear singular integral operator
with a 
Calder\'on-Zygmund kernel.
Then, $T$ can be extended to a bounded operator on $L^2$
if and only if

\begin{itemize}
\item[\textup{(i)}]
$ T $ satisfies the weak boundedness property,

\smallskip
\item[\textup{(ii)}]
$ T(1)$
and $T^*(1)$
are in $\BMO$.
\end{itemize}

\end{oldtheorem}

\noi
Since $T$ is a priori defined only in $\S$,
the expressions $T(1)$ and $T^*(1)$
are, of course, not  well defined and need to be interpreted carefully.
The same comment applies to the corresponding theorems in the bilinear setting.

The main concept needed in extending the $T(1)$ theorem to the $T(b)$ theorem is that of para-accretive functions.

\begin{definition}\label{DEF:accretive}\rm
We say that a function $b \in L^\infty$ is {\it para-accretive}\footnote{
An extra condition that $b^{-1}\in L^\infty$ is sometimes included in  the definition
of para-accretivity.
This, however, is not necessary.
Indeed, it follows from \eqref{AC2} and Lebesgue differentiation theorem
that $|b(x)| \geq c_0$ almost everywhere.
In particular, we have $b^{-1} \in L^\infty$.
}
if
there exists $c_0> 0$
such that, for every cube $Q$,
there exists a subcube $\wt Q \subset Q$
such that
\begin{align}
 \frac{1}{|Q|} \bigg|\int_{\wt Q} b(x) dx \bigg| \geq c_0.
\label{AC1}
 \end{align}

\end{definition}

\noi
It follows  from \eqref{AC1} that
\begin{align}
|\wt Q| \geq \frac{c_0}{\|b\|_{L^\infty}}|Q|.
\label{AC2}
\end{align}
In particular, the function 1 is automatically para-accretive. It is also worth pointing out that the definition of para-accretivity in the Definition \ref{DEF:accretive} is not the same as the one used in the classical $T(b)$ theorem of
David, Journ\'e, and Semmes \cite{DJS}. The notion of para-accretivity stated here is borrowed from
\cite{Han, Hart2}; for a similar definition in which cubes are replaced by balls, see Christ's monograph \cite{Christ}. The two definitions of para-accretivity are nevertheless equivalent. Since this natural observation seems to be missing from the literature, for the convenience of the reader, we have included its proof in the appendix.

Before giving a meaning to operators to which the $T(b)$ theorem applies, we need one more definition.

\begin{definition}\rm
Given $0< \eta \leq 1$, let $C^\eta$ be the collection of all functions from $\R^d\to \C$ such that $\|f\|_{C^\eta} < \infty$,
where the $C^\eta$-norm is given by
\[ \|f\|_{C^\eta}  = \sup_{x\ne y}\frac{|f(x) - f(y)|}{|x-y|^\eta}.\]
\noi
We also denote by $C^\eta_0$ the subspace of all compactly supported functions in $C^\eta$.
\end{definition}

\begin{definition}\label{DEF:CZ}\rm
Let $b_0$ and $b_1$ be para-accretive functions.
A linear singular operator $T:b_1C_0^\eta \to (b_0C_0^\eta)'$
is called  a \emph{linear singular integral operator
of Calder\'on-Zygmund type
associated to $b_0$ and $b_1$}
if $T$ is
continuous from
 $b_1C_0^\eta$ into $(b_0C_0^\eta)'$ for some $\eta>0$
 and there exists a \CZ kernel $K$
 such that
 \[ \jb {T (M_{b_1} f), b_0 g} = \int_{\R^{2d}} K(x, y) b_1(y) f(y) b_0(x) g(x) dx dy,\]

\noi
for all $f, g \in C_0^\eta$
such that $\supp f \cap \supp g = \emptyset$.
Here, $M_b$ denotes the operation of multiplication by $b$.

\end{definition}

With these preparations, we are now ready to state the classical $T(b)$ theorem \cite{DJS}.

\begin{oldtheorem}[$T(b)$ theorem]\label{THM:DJS}
Let $b_0$ and $b_1$ be para-accretive functions.
Suppose that $T$ is a linear singular integral operator
of Calder\'on-Zygmund type
associated to $b_0$ and $b_1$.
Then, $T$ can be extended to a bounded operator
on $L^2$
if and only if
the following conditions hold:
\begin{itemize}
\item[\textup{(i)}]
$M_{b_0} T M_{b_1}$ satisfies the weak boundedness property,

\smallskip
\item[\textup{(ii)}]
$ M_{b_0}T(b_1)$ and
$ M_{b_1}T^*(b_0)$
are in $\BMO$.

\end{itemize}

\end{oldtheorem}

\noi  In the special case when $b_0$ and $b_1$ are accretive\footnote{A function $b \in L^\infty$
is called accretive if there exists $\dl > 0$ such that $\Re b \geq \dl$ for all $x \in \R^d$.
Note that an accretive function is para-accretive.}
and $Tb_1 = T^*b_0 = 0$, the $T(b)$ theorem was independently proved by McIntosh and Meyer \cite{McM}.

\begin{remark}\label{REM:Hart}\rm
In \cite{DJS}, the condition (ii) of
Theorem \ref{THM:DJS}
is stated slightly differently; it was assumed that
$ T(b_1), T^*(b_0)\in \BMO$.
We note that this is just a matter of notation.
For example, the condition
$ T(b_1)\in \BMO$
in \cite{DJS}
means that  that there exists $\beta \in \BMO$ such that
\[\jb{ T(b_1), f} = \jb{\beta, f}
\quad
\text{for all mean-zero }f \in b_0 C_0^\eta.\]
This is clearly equivalent to
\begin{align}
\jb{ T(b_1), b_0 f} = \jb{\beta,b_0  f}
\quad
\text{for all }f \in C_0^\eta
\text{ such that } \int b_0 f dx = 0.
\label{Tb1}
\end{align}

\noi
Here, we used the fact that $b_0 f \leftrightarrow f$
is a one-to-one correspondence
since $b_0$ is para-accretive
and thus, in particular, is bounded away from zero almost everywhere.
In Theorem \ref{THM:DJS},
we followed the notation from \cite{Hart2}
to signify the fact that the condition indeed depends on both $b_0$ and $b_1$,
and what we mean by the condition (ii) in
 Theorem \ref{THM:DJS}
 is precisely the statement \eqref{Tb1}.
See also Theorem \ref{THM:BTB} below in the bilinear setting.

Lastly, note that, as in the $T(1)$ theorem,
the expressions
$ M_{b_0}T(b_1)$ and
$M_{b_1} T^*(b_0)$
are not a priori well defined
and thus some care must be taken.

\end{remark}

\subsection{Formulations of the $T(1)$ and $T(b)$  theorems \`a la Stein}

There is another formulation of the $T(1)$ theorem due to Stein
\cite{Stein} in which the conditions (i) and (ii) in Theorem \ref{THM:DJ}
are replaced by the quantitative estimate \eqref{S1} involving normalized bump functions.

\begin{oldtheorem}[$T(1)$ theorem \`a la Stein]\label{THM:T1Stein}
Let $T$ be as in Theorem \ref{THM:DJ}.
Then, $T$ can be extended to a bounded operator
on $L^2$
if and only if
there exists $M \in \mathbb{N} \cup\{0\}$
such that
we have
\begin{align}
\label{S1}
  \|T(\phi^{x_0, R})\|_{L^2}+\|T^*(\phi^{x_0, R})\|_{L^2}\les R^\frac{d}{2}
\end{align}

\noi
for any normalized bump function $\phi$
of order $M$, $x_0 \in \R^d$,
and $R>0$.
\end{oldtheorem}

By viewing the expressions
$T(\phi^{x_0, R})$
and
$T^*(\phi^{x_0, R})$
as
$T(1\cdot \phi^{x_0, R})$
and
$T^*(1\cdot \phi^{x_0, R})$,
it is  natural to extend
this result by replacing the constant function 1
by
para-accretive functions $b_0$ and $b_1$.
This is the first result of our paper.

\begin{maintheorem}[$T(b)$ theorem \`a la Stein]\label{THM:TbStein}
Let $T$, $b_0$, and  $b_1$ be as in Theorem \ref{THM:DJS}.
Then, $T$ can be extended to a bounded operator
on $L^2$
if and only if
there exists $M \in \mathbb{N} \cup\{0\}$
such that
the following two inequalities hold for any normalized bump function $\phi$
of order $M$, $x_0 \in \R^d$,
and $R>0$:
\begin{align}
  \|T(b_1\phi^{x_0, R})\|_{L^2}\les R^\frac{d}{2},  \label{TbS1}\\
   \|T^*(b_0\phi^{x_0, R})\|_{L^2}\les R^\frac{d}{2} \label{TbS2}.
\end{align}

\end{maintheorem}
	
\noi
We present the proof of Theorem \ref{THM:TbStein} in Section \ref{SEC:proof1}. 

As an application of this result, one could recover the well known fact that the commutator of a pseudodifferential operator with symbol in the H\"ormander class $S_{1, 0}^1$ and the multiplication operator of a Lipschitz function $a$ is bounded on $L^2$. Indeed, suppose that for all $x, \xi\in \R^d$ and all multi-indices $\alpha, \beta$ we have
\[|\partial_x^\alpha\partial_\xi^\beta \sigma (x, \xi)|\lesssim (1+|\xi|)^{1-|\beta|},\]
and let \[T_\sigma (f)(x)=\int_{\R^d} \sigma (x, \xi)f(\xi)e^{ix\cdot\xi}\,d\xi\] be the corresponding pseudodifferential operator. Also, given $a$ such that $\partial a/\partial x_j\in L^\infty (\R^d)$ for $1\leq j\leq d$, let
\[[T_\sigma, M_a]=T_\sigma (af)-aT(f)\] be the commutator of $T_\sigma$ and the multiplication operator $M_a$. It is straightforward to check that the kernel of $[T_\sigma, M_a]$ is Calder\'on-Zygmund and, by a similar computation to the one in \cite[pp.\,309-310]{Stein}, \eqref{TbS1} and \eqref{TbS2} hold as well; thus proving $[T_\sigma, M_a]: L^2\to L^2$.

\section{Bilinear \CZ theory}

Next, we turn our attention to the bilinear setting
and consider the corresponding extensions of the results in Section \ref{SEC:linear}.
Namely, we consider a bilinear singular operator
$T$ a priori defined from $\S\times\S$ into $ \S'$
of the form:
\begin{align}
 T(f, g)(x) = \int_{\R^{2d}} K(x, y, z) f(y) g(z) dy dz,
\label{BCZ1}
 \end{align}
\noi
where  we assume that, away from the diagonal
 $\Dl = \{ (x, y, z) \in \R^{3d}:\, x = y = z\}$,
 the distributional kernel $K$ coincides with
a function that is locally integrable on $\R^{3d} \setminus \Dl$.
The formal transposes
 $T^{*1}$ and $T^{*2}$
are defined in an analogous manner with
the
 kernels $K^{*1}$
and $K^{*2}$
given by \(K^{*1}(x, y, z) : = K(y, x, z)\,\, \text{and}\,\,
K^{*2}(x, y, z) : = K(z, y, x).\)

\begin{definition} \label{DEF:B}
\rm

A locally integrable function $K$ on $ \R^{3d} \setminus \Dl$
is called a  (bilinear) {\it   \CZ  kernel}
if it satisfies the following conditions.
\begin{itemize}
\item[\textup{(i)}]
For all $x, y, z \in \R^d$, we have
\[|K(x, y, z)|  \les \big(|x-y|+ |x-z|\big)^{-2d}, \]

\smallskip

\item[\textup{(ii)}]
There exists $\dl \in (0, 1]$ such that
\begin{align}
|K(x, y, z) - K(x', y, z)|
\les \frac{|x - x'|^\dl}{\big(|x-y|+ |x-z|\big)^{2d+\dl}}
\label{BCZ2}
\end{align}

\noi
for all $x, x', y,  z\in \R^d$
satisfying
$|x - x'| < \frac 12 \max\big(|x-y|, |x-z|\big)$.	
Moreover, we assume that
the formal transpose kernels $K^{*1}$
and $K^{*2}$
also
satisfy the regularity condition \eqref{BCZ2}.

\end{itemize}

\end{definition}

\noi
We say that a bilinear singular operator $T$
of the form \eqref{BCZ1}
with a bilinear \CZ kernel
is a {\it bilinear \CZ operator} if $T$ extends
to a bounded operator on $L^{p_0}\times L^{q_0}$ into $L^{r_0}$
for some $1< p_0 , q_0 < \infty$ with $\frac{1}{p_0} + \frac{1}{q_0} = \frac{1}{r_0} \leq 1$.

Similarly to the linear case, the crux of the bilinear \CZ theory is contained in the fact that if $T$ is a bilinear \CZ operator, then it is bounded on  $L^{p}\times L^{q}$ into $L^{r}$ for all $1< p , q < \infty$ with $\frac{1}{p} + \frac{1}{q} = \frac{1}{r} \leq 1$ (with the appropriate statements at the endpoints); see Grafakos and Torres \cite{GT1}. Therefore, the main question is to prove that there exists at least one triple $(p_0, q_0, r_0$) with  $1< p_0 , q_0 < \infty$ and $\frac{1}{p_0} + \frac{1}{q_0} = \frac{1}{r_0} \leq 1$
such that $T$ is bounded from  $L^{p_0}\times L^{q_0}$ into $L^{r_0}$.

The weak boundedness property for bilinear singular operators has a similar flavor as the one in the linear case.

\begin{definition} \label{DEF:WBP2} \rm

We say that a bilinear singular integral operator $T: \S\times \S\to \S'$
has the (bilinear) {\it weak boundedness property}
if there exists $M \in \mathbb N\cup \{0\}$
such that we have
\begin{align}
 \big| \jb{T(\phi_1^{x_1, R}, \phi_2^{x_2, R}), \phi_3^{x_3, R} } \big|\les R^d
\label{BWBP1}
 \end{align}

\noi
for all normalized bump functions $\phi_1, \phi_2, \phi_3$
of order $M$, $x_1, x_2, x_3 \in \R^d$,
and $R>0$.

\end{definition}

\begin{remark} \rm
It follows from  \cite[Lemma 9]{BO}
that it suffices to verify \eqref{BWBP1} for $x_1 = x_2 = x_3$.
\end{remark}

\subsection{Bilinear $T(1)$ and $T(b)$ theorems}
\label{SEC:bilinear}

We now state the bilinear $T(1)$ theorem in the form given by Hart \cite{Hart1}.

\begin{oldtheorem}[Bilinear $T(1)$ theorem]\label{THM:BT1}
Let $T: \S \times \S \to \S'$ be a bilinear singular integral operator
with a standard Calder\'on-Zygmund kernel.
Then, $T$ can be extended to a bounded operator on $L^p \times L^q  \to L^r$
for all $1 < p, q < \infty$ with $\frac 1p + \frac 1q = \frac 1r$
if and only if

\begin{itemize}
\item[\textup{(i)}]
$ T $ satisfies the  weak boundedness property,

\smallskip
\item[\textup{(ii)}]
$ T(1, 1), \,
T^{*1}(1, 1), $
and $ T^{*2}(1, 1)$
are in
$\BMO$.
\end{itemize}
\end{oldtheorem}

\noi
We chose  this formulation since it closely follows
the statement of the classical linear $T(1)$ theorem given in the previous section. Further, note that Theorem \ref{THM:BT1} is equivalent to the formulation of Grafakos-Torres \cite{GT1}; see also Christ and Journ\'e \cite{CJ}.

Next, we turn our attention to the bilinear version of the $T(b)$ theorem.

\begin{definition}\label{DEF:BCZ}\rm
 Let $b_0$, $b_1$, and $b_2$  be para-accretive functions.
A bilinear  singular operator $T:b_1C_0^\eta \times b_2C^\eta_0\to (b_0C_0^\eta)'$
is called   a \emph{bilinear singular integral operator
of Calder\'on-Zygmund type
associated to $b_0$,  $b_1$, and $b_2$}
if $T$ is
continuous
 from $b_1 C_0^\eta \times b_2 C_0^\eta $ into $(b_0 C_0^\eta)'$
 for some $\eta > 0$
  and there exists a bilinear \CZ kernel $K$
 such that
 \begin{equation}
  \jb {T (M_{b_1} f_1), M_{b_2} f_2), b_0 f_0}
 = \int_{\R^{3d}} K(x, y, z) b_0(x) f_0(x) b_1(y) f_1(y) b_2(z) f_2(z)  dx dy dz,
\label{tri1}
 \end{equation}

\noi
for all $f_0, f_1, f_2 \in C_0^\eta$ with
$\supp f_0\cap \supp f_1\cap \supp f_2= \emptyset$.

\end{definition}

Hart \cite{Hart2} proved the following result.

\begin{oldtheorem}[Bilinear $T(b)$ theorem]\label{THM:BTB}
 Let $b_0, b_1$,  and $b_2$ be para-accretive functions.
Suppose that $T$ is a bilinear singular integral operator
of Calder\'on-Zygmund type
associated to $b_0, b_1$, and $b_2$.
Then, $T$ can be extended to a bounded operator
on $L^p \times L^q \to L^r$
for all $p, q < \infty$ satisfying $\frac1p + \frac 1q = \frac 1r$
if and only if
\begin{itemize}
\item[\textup{(i)}]
$M_{b_0} T \big(M_{b_1}(\cdot), M_{b_2}(\cdot)\big) $ satisfies the weak boundedness property,

\smallskip
\item[\textup{(ii)}]
$M_{b_0}  T(b_1, b_2), \,
M_{b_1} T^{*1}(b_0, b_2)$,
and
$M_{b_2} T^{*2}(b_1, b_0)$
are in $\BMO$.

\end{itemize}

\end{oldtheorem}

\noi
As in Theorem \ref{THM:DJS},
we used the notation such as $M_{b_0}  T(b_1, b_2) \in \BMO$
rather than $  T(b_1, b_2)\in \BMO$
to signify the fact that each of the three statements in the condition (ii) of Theorem \ref{THM:BTB}
involves $b_0, b_1$ and $b_2$.
See Remark \ref{REM:Hart}.

\subsection{Formulation of the bilinear $T(b)$ theorem \`a la Stein}

As in the linear setting, we consider the formulation after Stein (Theorem \ref{THM:T1Stein}),
involving quantitative estimates.
In the following, we only state and prove  the formulation after Stein
in the context the bilinear $T(b)$ theorem.
The corresponding version for the bilinear $T(1)$ theorem
follows by setting $b_0 = b_1 = b_2 = 1$; this result already appears in \cite{GT1}.

\begin{maintheorem}[Bilinear $T(b)$ theorem \`a la Stein]\label{THM:BTbStein}
 Let $b_0, b_1$,  and $b_2$ be para-accretive functions.
Suppose that $T$ is a bilinear singular integral operator
of Calder\'on-Zygmund type
associated to $b_0, b_1$,  and $b_2$.
Then, $T$ can be extended to a bounded operator
on $L^p \times L^q \to L^r$
for all $p, q < \infty$ satisfying $\frac1p + \frac 1q = \frac 1r$
if and only if
there exists $M \in \mathbb{N} \cup\{0\}$
such that
we have
\begin{align}
  \|T(b_1\phi^{x_1, R}, b_2\phi^{x_2, R})\|_{L^2}\les R^\frac{d}{2},  \label{BTbS1}\\
  \|T^{*1}(b_0\phi^{x_0, R}, b_2\phi^{x_2, R})\|_{L^2}\les R^\frac{d}{2},  \label{BTbS2}\\
  \|T^{*2}(b_1\phi^{x_1, R}, b_0\phi^{x_0, R})\|_{L^2}\les R^\frac{d}{2}.  \label{BTbS3}
\end{align}

\noi
for any normalized bump function $\phi$
of order $M$, $x_0, x_1, x_2 \in \R^d$,
and $R>0$.

\end{maintheorem}

We prove this result in Section~\ref{proof2}.

\section{Proof  of Theorem \ref{THM:TbStein}}
\label{SEC:proof1}

Suppose that $T$ is bounded on $L^2$.
Let $\phi$ be a normalized bump function.
Then,
given any $x_0 \in \R^d$ and $R>0$,
we have
\begin{align*}
 \|T(b_1\phi^{x_0, R})\|_{L^2}
 \les
  \|b_1\|_{L^\infty}\|\phi^{x_0, R}\|_{L^2}
  \les R^\frac{d}{2}.
\end{align*}

\noi
This proves \eqref{TbS1}.
The condition \eqref{TbS2} follows from a similar computation.

Next, we assume that the conditions \eqref{TbS1} and \eqref{TbS2} hold. It suffices to show that the conditions \eqref{TbS1} and \eqref{TbS2} imply
the conditions (i) and  (ii) in Theorem \ref{THM:DJS}.

We first prove the condition (i) in Theorem \ref{THM:DJS}.
Let $\phi_1$ and $\phi_2$ be  normalized bump functions  of order 0.
Then, it follows from \eqref{TbS1} and \eqref{scaling} that we have
\begin{align*}
 \big| \big\langle M_{b_0}T M_{b_1}( \phi_1^{x_1, R}), \phi_2^{x_2, R}\big\rangle \big|
 \les \|b_0\|_{L^\infty} \|T (b_1 \phi_1^{x_1, R})\|_{L^2}\|\phi_2^{x_2, R}\|_{L^2}
 \les R^d
 \end{align*}

\noi
for
all  $x_1, x_2 \in \R^d$
and $R>0$.
This proves the weak boundedness property of $M_{b_0}TM_{b_1}$.

Next, we prove the condition (ii) in Theorem \ref{THM:DJS}.
In the following, we only show $M_{b_0}T(b_1) \in \BMO$,
assuming \eqref{TbS1}.
The proof of  $M_{b_1} T^*(b_0) \in \BMO$
follows from \eqref{TbS2} in an analogous manner.

We first recall from \cite{DJS} how to extend the definition of $T$
to $ b_1 C^\eta_b$, where $C^\eta_b := C^\eta\cap L^\infty$.
Denote by
$\{b_0C^\eta_0\}_0$
the subspace of mean-zero functions in $b_0C^\eta_0$.
Given $ f\in b_1 C^\eta_b$
and $g \in \{b_0C^\eta_0\}_0$,
let $\psi \in C^\eta_0$
with $0 \leq \psi \leq 1$
and $\psi \equiv 1$ in a neighborhood of $\supp  g$.
Then, we define the action of $T(f)$ on $g$ by
\begin{align}
\jb{T(f), g}
:& = \jb{T(f \psi ), g} + \jb{T(f(1-\psi)), g} \notag \\
& =  \jb{T(f \psi ), g} + \int_{\R^{2d}} \big[K(x, y) - K(x_0, y)\big]  f(y)(1-\psi(y)) g(x)dx dy.
\label{L1}
\end{align}

\noi
Note that this definition is independent of the choice of $\psi$.
Here, the last equality in \eqref{L1} holds
for any $x_0 \in \supp g$.

Let $\phi \in C^\infty_0$ with $0 \leq \phi \leq 1$
such that  $\phi (x) =  1$ for $|x| \leq \frac{1}{2}$ and $\supp  \phi \subset B_0(1)$.
Let $\phi_R(x) = \phi(R^{-1} x)$.
Then,
$T(b_1\phi_R)$ converges to $T(b_1)$ in the weak-$\ast$ topology
of $(\{b_0C^\eta_0\}_0)'$.
Namely, for all $g \in \{b_0C^\eta_0\}_0$, we have
\begin{align}
\jb{T(b_1), g} =
\lim_{R \to \infty}
\jb{T(b_1\phi_R), g}.
\label{L2}
\end{align}
Indeed,
letting $\psi \in C^\infty_0$ such that $\psi \equiv 1$ on $\supp g$ as before,
we have
\begin{align}
\jb{T(b_1\phi_R), g}
= \jb{T(b_1\psi \phi_R), g}
+ \jb{T(b_1(1-\psi) \phi_R), g}
\label{L3}
\end{align}

\noi
First, note that
\begin{align}
 \jb{T(b_1\psi \phi_R), g}
= \jb{T(b_1\psi), g}
\label{L4}
\end{align}

\noi	
for all sufficiently large $R$.
In view of \eqref{CZ2}, it follows from Lebesgue dominated convergence theorem
that
\begin{align}
\lim_{R \to \infty}
\jb{T& (b_1 (1-\psi)  \phi_R), g}
 =
 \lim_{R \to \infty}\int_{\R^{2d}} K(x, y) b_1(y) (1-\psi(y))  \phi_R(y)   g(x) dy dx \notag \\
& =
\lim_{R \to \infty}\int_{\R^{2d}} \big[K(x, y)- K(x_0, y)\big] b_1(y) (1-\psi(y))  \phi_R(y)  g(x) dy dx \notag \\
& = \int_{\R^{2d}} \big[K(x, y)- K(x_0, y)\big] b_1(y) (1-\psi(y))  g(x) dy dx \notag \\
& = \jb{T(b_1 (1-\psi)), g}
\label{L5}
\end{align}

\noi
where $x_0 \in \supp g$.
Then, \eqref{L2} follows from \eqref{L1}, \eqref{L3}, \eqref{L4}, and \eqref{L5}.

Suppose now that we have
\begin{align}
\|T(b_1 \phi_R)\|_{\BMO} \les 1,
\label{L6}
\end{align}

\noi
uniformly in $R>0$.
Then, by Banach-Alaoglu theorem
with $\BMO = (H^1)'$,
there exists a sequence $\{R_j\}_{j = 1}^\infty$
such that
$T(b_1 \phi_{R_j})$ converges in the weak-$\ast$ topology
to some function $\beta$ in $\BMO$.
Namely,
\begin{align}
\lim_{j \to \infty}
\jb{T(b_1\phi_{R_j}), g}
= \jb{\beta, g}
\label{L7}
\end{align}

\noi
for all $g \in H^1$.
In particular, \eqref{L7}  holds for all $g \in \{ b_0 C_0^\eta\}_0$.
Then, from \eqref{L2} and \eqref{L7}
with the uniqueness of a limit,
we can identify $T(b_1)$ (or rather $M_{b_0} T(b_1)$) with $\beta \in \BMO$. See Remark \ref{REM:Hart}.

Therefore, it remains to prove
\eqref{L6}.
Let $M \in \mathbb{N}\cup \{0\}$ be as in Theorem \ref{THM:TbStein}.
Then, by imposing that $\|\dd^\al \phi\|_{L^\infty} \leq 1$
for all multi-indices $\al$ with $|\al| \leq M$,
the function $\phi$ defined above is a normalized bump function of order $M$.

Fix a  cube $Q$ of side length $\l  > 0$ with center $x_0 \in\R^d$.
Set  $\phi_Q: = \phi^{x_0, r}$,
where
$r := 6 \diam(Q) = 6 \sqrt d \,  \l $.
By
writing $T(b_1 \phi_R)$ as
\begin{align}
 T(b_1\phi_R) = T(b_1 \phi_Q\phi_R) + T\big(b_1 (1-\phi_Q)\phi_R\big),
\label{L8}
 \end{align}

\noi
we consider the first and second terms separately.

On the one hand,
when $R \leq r$, write $\phi_Q\phi_R$ as
\begin{align}
\phi_Q(x)\phi_R(x)
=
  \phi(\tfrac{R}{r}\tfrac{x}{R}-\tfrac{x_0}{r}\big)\phi\big(\tfrac{x}{R}\big)
= [  \psi_1 \phi]^{0, R}(x)
\label{L8a}
\end{align}

\noi
with $\psi_1 (x) := \phi(\tfrac{R}{r}x-\tfrac{x_0}{r}\big)$.
Note that $ \psi_1 \phi$
is a normalized bump function.
Then, by the Cauchy-Schwarz inequality and \eqref{TbS1}, we have
\begin{align}
\int_{Q} \big|T(b_1 \phi_Q \phi_R) \big| dx
\leq  |Q|^\frac{1}{2} \big\|T\big(b_1 [\psi_1 \phi]^{0, R}\big)\big\|_{L^2}
\les R^\frac{d}{2} |Q|^\frac{1}{2}
\les |Q|.
\label{L9a}
\end{align}
	
\noi	
On the other hand,
when $R > r$, write $\phi_Q\phi_R$ as
\begin{align}
\phi_Q(x) \phi_R(x)
=
  \phi(\tfrac{x-x_0}{r}\big) \phi\big(\tfrac{r}{R} \tfrac{x-x_0}{r} + \tfrac{x_0}{R}\big)
= [\phi\psi_2]^{x_0, r}(x)
\label{L8b}
\end{align}

\noi
with $\psi_2 (x) :=
 \phi\big(\tfrac{r}{R} x + \tfrac{x_0}{R}\big)$.
Then, noting  that $\phi \psi_2$
is a normalized bump function,
it follows from the Cauchy-Schwarz inequality and \eqref{TbS1} that
\begin{align}
\int_{Q} \big|T(b_1 \phi_Q \phi_R) \big| dx
\leq  |Q|^\frac{1}{2} \big\|T\big(b_1 [\phi \psi_2]^{x_0, r}\big)\big\|_{L^2}
\les R^\frac{d}{2} |Q|^\frac{1}{2}
\les |Q|.
\label{L9b}
\end{align}

Next, we estimate the second term
in \eqref{L8}.
From the support condition:
\begin{align}
\supp (1- \phi_Q) \subset \R^d \setminus B_{x_0}( 3  \diam(Q))
\subset \R^d \setminus Q,
\label{L10}
\end{align}

\noi
we have
\[ T\big(b_1 (1-\phi_Q)\phi_R\big)(x) = 	
\int_{\R^d}  K(x, y) b_1(y) \big(1- \phi_Q(y) \big)
\phi_R(y)
dy,\]

\noi
for all $x\in Q$.
 Define $c_{Q, R}$ by
\[ c_{Q, R} : = \int_{\R^d}  K(x_0, y) b_1(y) \big(1- \phi_Q(y)\big)\phi_R(y) dy,\]

\noi
where  $x_0$ is the center of the cube $Q$.
Then,
it follows from \eqref{CZ2} with \eqref{L10} that,
 for $x \in Q$, we have
\begin{align}
\big| T\big(b_1 (1-\phi_Q) \phi_R\big)(x) - c_{Q, R}\big|
& \leq \int_{|x - x_0| \leq \diam(Q) \leq \frac 12 |x - y|}  |K(x, y) - K(x_0, y)| dy \notag \\
& \les 1
\label{L11}
\end{align}

\noi
uniformly in $R> 0$.

Hence,
putting  \eqref{L8}, \eqref{L9a}, \eqref{L9b},  and \eqref{L11} together,
we conclude that  there exists $A> 0$ such that for each cube $Q$
and $R> 0$,
there exists a constant $c_{Q, R}$ such that
\begin{align}
 \frac{1}{|Q|}\int_Q | T(b_1\phi_R)(x) - c_{Q, R}\big| dx \leq A.
\end{align}

\noi
Therefore,  it follows from
Proposition 7.1.2  in \cite{G2}
that
\[ \sup_{R>0} \|T(b_1\phi_R) \|_{\BMO} \leq 2A.\]

\noi
This proves \eqref{L7}
and thus completes the proof of Theorem \ref{THM:TbStein}.

\section{Proof of Theorem \ref{THM:BTbStein}}\label{proof2}

Suppose that $T$ is bounded on $L^4\times L^4 \to L^2$.
Then, given  a normalized bump function $\phi$,
we have
\begin{align*}
 \|T(b_1\phi^{x_1, R}, b_2\phi^{x_2, R} )\|_{L^2}
 \les
  \|b_1\|_{L^\infty}\|\phi^{x_1, R}\|_{L^4}
  \|b_2\|_{L^\infty}\|\phi^{x_2, R}\|_{L^4}
  \les R^\frac{d}{2}
\end{align*}

\noi
for any $x_1, x_2 \in \R^d$ and $R>0$.
This proves \eqref{BTbS1}.
A similar computation
yields \eqref{BTbS2} and \eqref{BTbS3}.

Next, we assume that the conditions
\eqref{BTbS1}, \eqref{BTbS2}, and \eqref{BTbS3} hold. It suffices to show that
the conditions \eqref{BTbS1}, \eqref{BTbS2},  and \eqref{BTbS3} imply
the conditions (i) and  (ii) in Theorem \ref{THM:BTB}.

We first prove the condition (i) in Theorem \ref{THM:BTB}.
Let $\phi_j$, $j = 0, 1, 2$,  be  normalized bump functions  of order $M$.
Then, it follows from the Cauchy-Schwarz inequality, \eqref{BTbS1},  and \eqref{scaling} that
\[
 \big| \jb{M_{b_0}T (b_1 \phi_1^{x_1, R}, b_2\phi_2^{x_2, R}),
  \phi_0^{x_0, R}} \big|  \les \|b_0\|_{L^\infty}
\| T (b_1 \phi_1^{x_1, R}, b_2\phi_2^{x_2, R})\|_{L^2}
\|  \phi_0^{x_0, R} \|_{L^2}
 \les R^d
\]

\noi
for
all  $x_0, x_1, x_2 \in \R^d$
and $R>0$.
This proves the condition (i) in Theorem \ref{THM:BTB}.

Next, we prove the condition (ii) in Theorem \ref{THM:BTB}.
As in the proof of Theorem \ref{THM:TbStein}, we only show
$M_{b_0}T(b_1, b_2) \in \BMO$,
assuming \eqref{BTbS1}.
The proof of the other two conditions
follows in a similar manner in view of the symmetric condition in
Definition \ref{DEF:B}.

Since $T$ is a priori defined only on
 $ b_1 C^\eta_0 \times b_2 C_0^\eta$,
 we first extend $T$ to  $ b_1 C^\eta_b \times b_2 C_b^\eta$.
Fix $ f_j\in b_j C^\eta_b$ , $j = 1, 2$.
Given  $g \in \{b_0C^\eta_0\}_0$,
let $\psi \in C^\eta_0$
with $0 \leq \psi \leq 1$
and $\psi \equiv 1$ in a neighborhood of $\supp  g$.
Then, we define the action of $T(f_1, f_2)$ on $g$ by
\begin{align}
\jb{T(f_1, f_2), g}
& := \jb{T(f_1 \psi, f_2 \psi ), g}
+
\jb{T(f_1 (1-\psi), f_2 \psi ), g} \notag \\
& \hphantom{XX}
+ \jb{T(f_1 \psi, f_2 (1-\psi) ), g}
+ \jb{T(f_1 (1-\psi), f_2 (1-\psi) ), g}.
\label{B1}
\end{align}

\noi
Note that  the last three terms can be written as triple integrals of the form \eqref{tri1}.
From this, we see that this definition is independent of the choice of $\psi$.

Let $\phi \in C^\infty_0$ with $0 \leq \phi \leq 1$
such that  $\phi (x) =  1$ for $|x| \leq \frac{1}{2}$ and $\supp  \phi \subset B_0(1)$.
Let $\phi_R(x) = \phi(R^{-1} x)$.
Then,
$T(b_1\phi_R, b_2\phi_R)$ converges to $T(b_1, b_2)$ in the weak-$\ast$ topology
of $(\{b_0C^\eta_0\}_0)'$.
Namely, we have
\begin{align}
\jb{T(b_1, b_2), g} =
\lim_{R \to \infty}
\jb{T(b_1\phi_R, b_2\phi_R), g}
\label{B2}
\end{align}

\noi
for all $g \in \{b_0C^\eta_0\}_0$.
See \cite{Hart2} for the proof of \eqref{B2}.

Suppose that we have
\begin{align}
\|T(b_1 \phi_R, b_2 \phi_R)\|_{\BMO} \les 1,
\label{B3}
\end{align}

\noi
uniformly in $R>0$.
Then, as in the proof of Theorem \ref{THM:TbStein},
it follows from Banach-Alaoglu theorem
that
there exists a sequence $\{R_j\}_{j = 1}^\infty$
and $\be \in \BMO$
such that
\begin{align}
\lim_{j \to \infty}
\jb{T(b_1\phi_{R_j}, b_2 \phi_{R_j}), g}
= \jb{\beta, g}
\label{B4}
\end{align}

\noi
for all $g \in H^1$, in particular for
 all $g \in \{ b_0 C_0^\eta\}_0$.
Hence, from \eqref{B2} and \eqref{B4},
we conclude that  $M_{b_0} T(b_1, b_2) \in \BMO$.

Therefore, it remains to prove
\eqref{B3}.
By imposing that $\|\dd^\al \phi\|_{L^\infty} \leq 1$
for all multi-indices $\al$ with $|\al| \leq M$,
the function $\phi$ defined above is
 a normalized bump function of order $M$.
As in the proof of Theorem \ref{THM:TbStein},
let $Q$ be the cube of side length $\l  > 0$ with center $x_0 \in\R^d$. Set  $\phi_Q = \phi^{x_0, r}$, where $r = 6 \diam(Q)$.
Then,
write $T(b_1 \phi_R, b_2 \phi_R)$ as
\begin{align}
 T(b_1\phi_R, b_2 \phi_R)
 & = T(b_1 \phi_Q\phi_R, b_2 \phi_Q\phi_R)
  + T\big(b_1 (1-\phi_Q)\phi_R, b_2 \phi_Q\phi_R\big)\notag \\
&  \hphantom{X}
 + T\big(b_1 \phi_Q\phi_R, b_2 (1-\phi_Q)\phi_R\big)
 + T\big(b_1 (1-\phi_Q)\phi_R, b_2 (1-\phi_Q)\phi_R\big)
 \notag\\
 & :=  \I+\II + \III + \IV.
\label{B5}
 \end{align}

It follows from the Cauchy-Schwarz inequality and \eqref{BTbS1}
with \eqref{L8a} and \eqref{L8b} that
\begin{align}
\int_{Q} |\I|  dx
\leq
\begin{cases}
\vphantom{\Big|} |Q|^\frac{1}{2} \big\|T\big(b_1 [\psi_1\phi]^{0, R}, b_2 [\psi_1\phi]^{0, R}\big)\big\|_{L^2}
\les |Q|,
 & \text{when }R\leq r, \\
\vphantom{\Big|}
 |Q|^\frac{1}{2} \big\|T\big(b_1 [\phi\psi_2]^{x_0, r}, b_2 [\phi\psi_2]^{x_0, r}\big)\big\|_{L^2}
\les |Q|,
 & \text{when }R> r.\\
\end{cases}
\label{B6}
\end{align}

Next, we consider the terms $\II, \III$, and $\IV$.
Let $\phi_Q^c := 1 - \phi_Q$.
Then,
from the support condition \eqref{L10},
we have
\[ \II(x) =
\int_{\R^{2d}}  K(x, y, z)
 b_1(y)  \phi^c_Q(y)\phi_R(y)
b_2(z) \phi_Q(z)
\phi_R(z)
dydz
\]

\noi
 for $x\in Q$.
Define $c^{(2)}_{Q, R}$ by
\[ c^{(2)}_{Q, R} : = \int_{\R^{2d}}  K(x_0, y, z)
 b_1(y)  \phi^c_Q(y)\phi_R(y)
b_2(z) \phi_Q(z)
\phi_R(z)
dydz, \]

\noi
where $x_0$ is the center of the cube $Q$.
Then, it follows from \eqref{BCZ2} with \eqref{L10} that, for $x \in Q$, we have
\begin{align}
| \II(x) - c^{(2)}_{Q, R}\big|
& \leq \int_{\supp \phi_Q}\int_{|x - x_0| \leq \diam(Q) \leq \frac 12 |x - y|}
  |K(x, y, z) - K(x_0, y, z)| dy dz \notag \\
& \les 1
\label{B7}
\end{align}

\noi
uniformly in $R> 0$.
By symmetry, the same estimate holds
for $\III$.
As for $\IV$,
by letting
\[ c^{(4)}_{Q, R} : = \int_{\R^{2d}}  K(x_0, y, z)
 b_1(y)  \phi^c_Q(y)\phi_R(y)
b_2(z) \phi_Q^c(z)
\phi_R(z)
dydz, \]

\noi
we have, for $x \in Q$,
\begin{align}
| \IV(x) - c^{(4)}_{Q, R}\big|
& \leq \int_{|x - x_0| \leq \diam(Q) \leq \frac 12 \min( |x - y|, |x-z|)}
  |K(x, y, z) - K(x_0, y, z)| dy dz \notag \\
& \les 1
\label{B8}
\end{align}

\noi
uniformly in $R> 0$.

Hence,
putting  \eqref{B5}, \eqref{B6}, \eqref{B7},  and \eqref{B8} together,
we conclude that  there exists $A> 0$ such that for each cube $Q$
and $R> 0$,
there exists a constant $\wt c_{Q, R}$ such that
\begin{align}
 \frac{1}{|Q|}\int_Q | T(b_1\phi_R, b_2\phi_R)(x) - \wt c_{Q, R}\big| dx \leq A,
\end{align}

\noi
thus yielding
 \eqref{B3}.
This completes the proof of Theorem \ref{THM:BTbStein}.

\appendix

\section{On para-accretive functions}

Para-accretive functions play an important role
in the $T(b)$ theorems.
In this paper, we used Definition \ref{DEF:accretive}
for para-accretivity.
In \cite{DJS}, however, David, Journ\'e, and Semmes
used a different definition (see Definition \ref{DEF:PA} below) and
gave several equivalent characterizations for para-accretive functions
(Proposition \ref{PROP:PA2} below). In this appendix, we
show that these two definitions (Definition \ref{DEF:accretive} and Definition \ref{DEF:PA}) are equivalent.

\begin{definition}\label{DEF:PA}\rm
A function $b \in L^\infty$ is para-accretive
if $b^{-1} \in L^\infty$
and there exists a sequence $\{s_k\}_{k\in \Z}$
of functions $s_k:\R^d \times \R^d \to \C$
for which the following conditions hold;
there exist $C>0$ and $\al > 0$ such that for all $k\in \Z$,
\begin{align*}
\text{(i)} & \quad  |s_k(x, y)| \leq C 2^{kd}, \quad \text{for all  }x, y \in \R^d,\\
\text{(ii)} &
\quad s_k(x, y) = 0,
\quad   \text{if  }
|x- y|\geq C2^{-k},\\
\text{(iii)} & \quad
s_k(x, y) = s_k(y, x),
\quad   \text{for all  } x, y \in \R^d, \\
\text{(iv)} & \quad
|s_k(x, y) -  s_k(x', y)|\leq C2^{k(d+\al)}|x-x'|,
\quad   \text{for all  } x, x', y \in \R^d, \\
\text{(v)} & \quad
\textstyle \int s_k(x, y)b(y) dy = 1,
\quad   \text{for all  } x \in \R^d.
\end{align*}

\end{definition}

\noi
The following proposition states
different  characterizations for para-accretive functions
 according to Definition \ref{DEF:PA}.

\begin{proposition}[Proposition 2 in \cite{DJS}]\label{PROP:PA2}
Let $b \in L^\infty$ such that $b^{-1} \in L^\infty$.
Then, the following statements are equivalent.

\smallskip
\noi
\textup{(A)}
A function $b$ is para-accretive according to Definition \ref{DEF:PA}.

\smallskip
\noi
\textup{(B)}
There exists $\eps>0$ and $N >0$
such that
for all $k \in \Z$
and for any  dyadic cube $Q$ of side length $\l(Q) = 2^{-k}$,
there exists another dyadic cube $\wt Q$ of the same side length
such that the distance between $Q$ and $\wt Q$ is at most $N 2^{-k}$
and
\[  \frac{1}{|\wt Q|} \bigg|\int_{\wt Q} b(x) dx\bigg| \geq \eps_1.\]

\smallskip
\noi
\textup{(C)}
There exist
$C >0 $, $\dl >0$,
and $u_k: \R^d\times \R^d \to \C$
such that for all $k \in \Z$,
\begin{align*}
\textup{(i)} & \quad  |u_k(x, y)| \leq C 2^{kd}, \quad \text{for all  }x, y \in \R^d,\\
\textup{(ii)} &
\quad u_k(x, y) = 0,
\quad   \text{if  }
|x- y|\geq C2^{-k},\\
\textup{(iii)} & \quad
|u_k(x, y)  - u_k(x, y')|\leq C2^{k(d+\dl)}|y-y'|,
\quad   \text{for all  } x, y, y' \in \R^d, \\
\textup{(iv)} &
\quad \text{For all  } x \in \R^d, \\
& \hphantom{XXXXX}
 \frac{1}{C} \leq \bigg|\int u_k(x, y)b(y) dy\bigg|\leq C.
\quad
\end{align*}

\smallskip
\noi
\textup{(D)}
There exist
$C >0 $, $\dl >0$,
and $v_k: \R^d\times \R^d \to \C$
such that for all $k \in \Z$,
the conditions \textup{(i)}-\textup{(iv)} in \textup{(C)} are satisfied.
Moreover,
the following extra conditions are satisfied:
\begin{align*}
\textup{(v)} & \quad
\int v_k(x, y) dy = 1,
\quad \text{for all  } x \in \R^d,\\
\textup{(vi)} &
\quad
\int v_k(x, y) dx = 1,
\quad \text{for all  } y \in \R^d,
\end{align*}

\noi
and
\begin{itemize}
\item[\textup{(vii)}]
For all $y \in \R^d$, the function $v_k(\,\cdot\,, y)$ is constant for each
dyadic cube of side length $2^{-k}$.
\end{itemize}

\end{proposition}


\smallskip

\subsection{Definition	
 \ref{DEF:PA}
implies Definition
 \ref{DEF:accretive}} Let   $b\in L^\infty$ be para-accretive
according to
Definition
 \ref{DEF:PA}.
 In the following, we show
that \eqref{AC1} in Definition
 \ref{DEF:accretive}
follows from   Proposition \ref{PROP:PA2} (B).

Let $\eps>0$ and $N>0$ be as in
 Proposition \ref{PROP:PA2} (B).
 Without loss of generality, we assume that $N\geq 10$.
Given  a dyadic cube $Q$ centered at $x_0$,
choose $k \in \Z$ such that
\begin{align}
10\cdot N2^{-k} \leq \l (Q) \leq 20\cdot N2^{-k}.
\label{M1}
\end{align}

\noi
Fix  a dyadic cube $Q_1\subset Q$  of side length $2^{-k}$,
containing the center $x_0$ of the cube $Q$.
Then, by  Proposition \ref{PROP:PA2} (B),
there exists another dyadic cube $Q_2$
of side length $2^{-k}$
within distance $N2^{-k}$ from $Q_1$
such that
\begin{align}
 \frac{1}{|Q_2|}\bigg|\int_{Q_2}b(x)dx\bigg|\geq \eps.
\label{M2}
 \end{align}

\noi
Note that
$Q_2\subset Q$.
Moreover, from \eqref{M1} and \eqref{M2}, we have
\begin{align*}
 \frac{1}{|Q|}\bigg|\int_{Q_2}b(x)dx\bigg|\geq \frac{\eps}{(20N)^d}.
 \end{align*}

\noi
Since the choice of $Q$ was arbitrary,
this shows that $b$ is indeed para-accretive in the sense
of Definition \ref{DEF:accretive}.

\medskip

\subsection{Definition	
 \ref{DEF:accretive}
implies Definition
 \ref{DEF:PA}} Let   $b\in L^\infty$ be para-accretive
according to
Definition
 \ref{DEF:accretive}.
It suffices to construct a sequence $\{u_k\}_{k\in \Z}$
of functions $u_k$ on $\R^d \times \R^d$,
satisfying the conditions (i)-(iv) in
 Proposition \ref{PROP:PA2} (C).

Let $\phi\in C_0^\infty$
be a normalized bump function of order 1
such that
$\int_{\R^d} \phi(x) dx = \al^{-1}>0$. Then, let $\phi_\eps(x) = \eps^{-d} \al \phi(\eps^{-1} x)$, that is, $\{ \phi_\eps\}_{\eps > 0}$ is an approximation to the identity.

Given $k \in \Z$,
let $Q_k$ be
the cube
 of side length $2^{-k}$
 centered at the origin
 and $Q^x_k  := x+ Q_k$
be the cube of side length $2^{-k}$
centered at $x \in \R^d$.
Then, by Definition \ref{DEF:accretive},
there exists a subcube $\wt Q^x_k \subset Q^x_k$
 such that
 \begin{align}
2^{kd} \bigg|\int \ind_{\wt Q^x_k}(y) b(y) dy \bigg| \geq c_0.
 \label{N1}
 \end{align}
	
\noi	
 Here, $c_0$ is uniform in all cubes $Q^x_k \supset \wt Q^x_k$.
 From \eqref{AC2}, we also have
 \begin{align}
\l_{x, k}:= \l (\wt Q^x_k) \geq c_1\l(Q^x_k)
= c_12^{-k},
\qquad \text{where }c_1 = c_1(b):= \bigg(\frac{c_0}{\|b\|_{L^\infty}}\bigg)^\frac{1}{d}.
 \label{N2}
 \end{align}

\noi

Note that
\begin{align}
\ind_{\wt Q^x_k}(y) = \ind_{Q_0}(\l_{x, k}^{-1} (y-\wt x)),
\label{N3}
\end{align}
	
\noi
where $\wt x$ is the center of the subcube $\wt Q^x_k$.
Then, by setting $\eps = h \l_{x, k}$ for $h > 0$,
we have
\begin{align}
\ind_{\wt Q^x_k}*\phi_\eps(y)
& = \int \ind_{\wt Q^x_k}(y-z)\phi_\eps(z)dz
= \l^{-d}_{x, k}
\int \ind_{Q_0}(\l^{-1}_{x, k}(y-\wt x-z))\phi_h(\l^{-1}_{x, k}z)dz \notag\\
& =
\ind_{Q_0}*\phi_h(\l^{-1}_{x, k}(y-\wt x)).
\label{N4}
\end{align}

\noi
Then,
it follows from from \eqref{N3}
and \eqref{N4}
that we can choose sufficiently small $h\ll1$
such that
\begin{align}
2^{kd} \bigg|\int
 \ind_{\wt Q^x_k}*\phi_\eps(y)
 - \ind_{\wt Q^x_k}(y)dy\bigg|
&  = \frac{|\wt Q^x_k|}{|Q^x_k|}
 \bigg|\int
\ind_{\wt Q_0}*\phi_h(y)
- \ind_{\wt Q_0}(y)dy\bigg|\notag \\
& \leq \bigg|\int\ind_{\wt Q_0}*\phi_h(y)
- \ind_{\wt Q_0}(y)dy\bigg|
\leq \frac{c_0}{2 \|b\|_{L^\infty}},
\label{N5}
\end{align}
	
\noi
uniformly in $x \in \R^d$ and $k \in \Z$.
Hence, using \eqref{N1}, \eqref{N5}, and the triangle inequality, we obtain
\begin{align}
2^{kd} \bigg|\int
 \ind_{\wt Q^x_k}*\phi_{h \l_{x, k}}(y)b(y) dy\bigg|
\geq \frac{1}{2}c_0.
\label{N6}
\end{align}

\noi
In the following, we fix
 $h \ll1$ such that \eqref{N5} holds.

Now, let us define $u_k$ by
\begin{align}
 u_k (x, y) := |Q_k|^{-1}  \ind_{\wt Q^x_k}*\phi_{h \l_{x, k}}(y)
= 2^{kd} \ind_{\wt Q^x_k}*\phi_{h \l_{x, k}}(y).
\label{N7}
 \end{align}

\noi
Then, from \eqref{N6} and Young's inequality, we have
\begin{align*}
\frac{1}{2}c_0  \leq \bigg|\int u_k(x, y)b(y) dy\bigg|
\leq \frac{|\wt Q^x_k|}{|Q^x_k|} \|\phi_{h\l_{x, k}}\|_{L^1} \|b\|_{L^\infty}
\leq  \|b\|_{L^\infty}
\end{align*}

\noi
for all $x\in \R^d$ and $k \in \Z$.
Hence, (iv) holds.

By the mean value theorem
and Young's inequality with \eqref{N2}, we have
\begin{align*}
|u_k(x, y) -u_k(x, y')|
& \leq
2^{kd}\|\ind_{\wt Q^x_k}*\dd(\phi_{h \l_{x, k}})\|_{L^\infty}
|y-y'|
 \leq \al 2^{kd}|\wt Q^x_k|
(h \l_{x, k})^{-d-1}
|y - y'| \notag \\
& \leq
\al c_1^{-1} h^{-d-1}
2^{k(d+1)}
|y - y'|
\end{align*}

\noi
for all $x, y, y'\in \R^d$.
This proves (iii).
By Young's inequality, we have
\begin{align*}
\|u_k(x, y)\|_{L^\infty}
\leq   2^{kd} \|\ind_{\wt Q^x_k}\|_{L^1} \|\phi_{h \l_{x, k}}\|_{L^\infty}
\leq \al 2^{kd}
|\wt Q^x_k|  (h\l_{x, k})^{-d} = \al h^{-d} 2^{kd}
\end{align*}

\noi
for all $x, y\in \R^d$.
This proves (i).
Lastly, from
\eqref{N4} and \eqref{N7}, we have
\begin{align}
 u_k (x, y)
= 2^{kd} \ind_{Q_0}*\phi_h(\l^{-1}_{x, k}(y-\wt x))
= 0
\label{N8}
 \end{align}

\noi
for $ |\wt x-y| \geq \big(1+ \tfrac{\sqrt{d}}{2}\big) \l_{x, k}$
since $h \leq 1$.
Note that
(i)
$\l_{x, k}= \l (\wt Q^x_k) \leq \l(Q^x_k) = 2^{-k}$
and (ii) $|x - \wt x| \leq \frac{\sqrt d}{2}2^{-k}$,
since $x$ and $\wt x$ are the centers of the cubes $Q_k^x$
and $\wt Q_k^x$, respectively.
This in particular implies
that \eqref{N8} holds
for $ | x-y| \geq \big(1+ \sqrt{d}\big) 2^{-k}$.
This proves the condition (ii).
By Proposition \ref{PROP:PA2},
we conclude that $b$ is para-accretive
in the sense of
Definition
 \ref{DEF:PA}.

Therefore,
Definitions
 \ref{DEF:accretive} and
 Definition
 \ref{DEF:PA}
 are equivalent.

\smallskip

\noi
{\bf Acknowledgments.}
This work was partially supported by a grant from the Simons Foundation (No.~246024 to \'Arp\'ad B\'enyi). The authors are grateful to
the Hausdorff Research Institute for Mathematics in Bonn
for its  hospitality during the trimester program Harmonic Analysis and Partial Differential Equations, where a part of this manuscript was prepared. They would also like to thank Jarod Hart for  helpful discussions.

\end{document}